\documentclass[12pt,a4paper,reqno]{article}
\usepackage{latexsym}
\usepackage{amssymb}
\usepackage{amsmath}
\usepackage{xspace}
\usepackage{enumerate}
\usepackage{amsthm}
\usepackage{fancyhdr}
\usepackage[all]{xy}
\newtheorem{thm}{Theorem}[section]
\newtheorem{prop}[thm]{Proposition}
\newtheorem{lem}[thm]{Lemma}
\theoremstyle{definition}
\newtheorem{cor}[thm]{Corollary}
\newtheorem{defn}[thm]{Definition}
\newtheorem*{nota}{Notation}

\newtheorem{exam}[thm]{Example}

\theoremstyle{remark}
\newtheorem*{rem}{Remark:}
\newtheorem*{rems}{Remarks:}

\newcommand{\coend}{\ensuremath{e_D(X)}\xspace}
\newcommand{\rmod}{\ensuremath{\textrm{$R$-Mod}}\xspace}
\newcommand{\tens}{\ensuremath{\otimes}}
\newcommand{\cotens}[1]{\ensuremath{\Box_{#1}}}
\newcommand{\cotensc}{\cotens{C}}
\newcommand{\cotensd}{\cotens{D}}
\newcommand{\comc}{\ensuremath{\textrm{Com}_C}}
\newcommand{\comd}{\ensuremath{\textrm{Com}_D}}
\newcommand{\homr}{\ensuremath{\textrm{Hom}_R}}
\newcommand{\cohomf}{\ensuremath{\mathit{h}_D(X,-)}}
\newcommand{\cohomx}[1]{\ensuremath{\mathit{h}_D(X,#1)}}
\newcommand{\cohomb}{\ensuremath{\mathit{h}(-,-)}}
\newcommand{\comod}[1]{\ensuremath{\mathcal{M}^#1}\xspace}
\newcommand{\comodc}{\comod{C}}
\newcommand{\comodd}{\comod{D}}
\newcommand{\ccomod}{\ensuremath{{}^C\mathcal{M}}\xspace}
\newcommand{\dcomod}{\ensuremath{{}^D\mathcal{M}}\xspace}
\newcommand{\cdbicomod}{\ensuremath{{}^C\mathcal{M}^D}}
\newcommand{\dcbicomod}{\ensuremath{{}^D\mathcal{M}^C}}
\newcommand{\bicomod}[2]{\ensuremath{{}^{#1}\mathcal{M}^{#2}}}
\newcommand{\rcomod}[1]{\ensuremath{\varrho_{#1}:#1\to #1\tens C}}
\newcommand{\lcomod}[1]{\ensuremath{\varrho_{#1}:#1\to C\tens #1}}
\newcommand{\coalg}[1]{\ensuremath{\Delta:#1\to #1\tens #1}}
\newcommand{\comodcf}{\ensuremath{\mathcal{M}^{C}_f}\xspace}
\newcommand{\comoddf}{\ensuremath{\mathcal{M}^{D}_f}\xspace}
\newcommand{\ccomodf}{\ensuremath{^{C}\mathcal{M}_f}\xspace}
\newcommand{\xra}{\xrightarrow}

\newcommand{\lra}{\longrightarrow}
\newcommand{\lam}[2]{\ensuremath{\lambda_{#1,#2}}}
\newcommand{\lamb}[2]{\ensuremath{\overline{\lambda}_{#1,#2}}}

\newcommand{\roh}{\varrho}

\begin{document}
\title{Equivalences of Comodule Categories for Coalgebras over Rings.}
\author{Khaled AL-Takhman
\\ \small{Department of Mathematics}
\\ \small{Birzeit University}
\\ \small{P. O. Box 14, Birzeit, Palestine} \\
  \small{E-mail: takhman@birzeit.edu}
}
\date{}
\maketitle
\begin{abstract}
In this article we defined and studied quasi-finite comodules, the cohom
functors for coalgebras over rings. linear functors between categories
of comodules are also investigated and it is proved that good enough
linear functors are nothing but a cotensor functor. Our main result of
this work characterizes equivalences between comodule categories
generalizing the Morita-Takeuchi theory to coalgebras over rings.
Morita-Takeuchi contexts in our setting is defined and investigated, a
correspondence between strict Morita-Takeuchi contexts and equivalences
of comodule categories over the involved coalgebras is obtained. Finally
we proved that for coalgebras over QF-rings Takeuchi's representation of
the cohom-functor is also valid.
\end{abstract}
\vspace{2cm} \emph{Keywords}:cotensor product, cohom functor,
faithfully coflat comodules, equivalent categories,
Morita-Takeuchi context.


\section*{Introduction}
M. Takeuchi \cite{Tak77} developed a theorem that characterizes
equivalences of comodule categories over fields, dualizing Morita
results on equivalences of module categories. In recent years a
new interest arose for the study of coalgebras over rings. In
this article we did this for equivalences of comodule categories.

One relevant consequence of the Morita theory on equivalences
between module categories is that the notion of Morita equivalent
rings is independent 
categories of right modules are). Here, we extend the
aforementioned equivalence theory to coalgebras over an arbitrary
commutative ring $R$, under the mild hypothesis that the involved
coalgebras are flat as $R$-modules. To transfer the
Morita-Takeuchi theory from fields to commutative rings, we
overcome two difficulties: on the one hand, no basis-dependent
arguments can be used here, and, on the other, the lack of
associativity in general of the cotensor product.

\vspace{.3cm} The cotensor functor plays the most important role
in characterizing equivalences of comodule categories. S.
Caenepeel \cite{Cae98} gave some sufficient conditions to
guarantee the associativity of the cotensor product over rings,
but unfortunately his result is not true (for a counter example
see \cite{Gru87}). A detailed study of the cotensor functor for
coalgebras over rings is done by the author in \cite{Altconf}. In
section \ref{s2} we state some results on the exactness and
associativity of this functor that we need in this work.

\vspace{.3cm}To prove our results we developed purely categorical
arguments for the study of quasi-finite comodules and  the
cohom-functor in section \ref{s3}. In this section we considered
the coendomorphism coalgebra derived from a quasi-finite comodule
$X$ and showed that the dual algebra of this coalgebra is algebra
(anti-) isomorphic to the algebra of comodule endomorphisms of
$X$. In section \ref{s4} we studied linear functors (in
particular equivalences) between comodule categories and showed
that a linear functor, under certain conditions, is isomorphic to
a cotensor functor. In section \ref{s5} we proved our main result
(Theorem \ref{mtt}) which characterizes an equivalence by the
existence of a bicomodule which is quasi-finite, faithfully
coflat and an injector on either side.

\vspace{.3cm}The concept of Morita-Takeuchi context was first
introduced in \cite{Tak77} for coalgebras over fields. In
\cite{Das95} it was studied for graded coalgebras. S. Caenepeel
\cite{Cae98} defined it for coalgebras over rings, in his
definition the purity condition which guarantees the
associativity of the cotensor product is overlooked, so his
definition for strict Morita-Takeuchi context does not give an
equivalence between the comodule categories of the involved
coalgebras. In section \ref{Scontext} we defined this concept for
coalgebras over rings, and showed that there is a correspondence
between strict Morita-Takeuchi contexts and equivalences of
comodule categories. We also showed that for a $C$-comodule $X$
that is quasi-finite, faithfully coflat and an injector the
coendomorphism coalgebra of $X$ is Morita-Takeuchi equivalent to
$C$.

\vspace{.3cm}The ground ring has great influence on the
properties of comodule categories. In section \ref{s6} we studied
the case where the ring is a QF-ring and proved that most of the
results of \cite{Tak77} are true for our settings.

\vspace{.3cm}Throughout this paper we assume that all rings are
commutative with unity, all modules are unitary, and the
unadorned tensor product is understood to be over the ground
ring. The categorical terminology we used are those of
\cite{Sten} with minor differences. For module theoretic notions
we refer to \cite{Wisbi}. Finally, a submodule $W$ of an
$R$-module $V$ is called $N$-pure ($N$ is an $R$-module) if the
canonical map $W\tens N\to V\tens N$ is a monomorphism. $W$ is
called pure if it is $N$-pure for every $N$.

\section{Notations and Preliminaries}
Let $R$ be a commutative ring with unity. We denote by $\rmod $
the category of unital $R$-modules. In this section we recall the 
basic definitions and results that we need in the sequel.
\par
\smallskip
\textbf{Coalgebras} : A {\em coalgebra} over a ring $R$ is an
$R$-module $C$ together with two $R$-linear maps $ \Delta : C\to
C\tens C $ ({\em comultiplication}) and $\varepsilon : C\to R $
({\em counit}) such that
 $ (id\tens \Delta)\circ \Delta = (\Delta \tens id)\circ \Delta
 $ ({\em Coassociativity property})  and $ (\varepsilon \tens id)
\circ\Delta = id = (id \tens \varepsilon)\circ \Delta $ ({\em
counitary property}).
\smallskip
\par
\textbf{Coalgebra Morphisms} : Let $ \Delta_C : C\to C\tens C $, $
 \Delta_D : D\to D\tens D $  be coalgebras, an $R$-linear map $\pi
: C \to D $ is called {\em a coalgebra morphism} if $
\Delta_D\circ \pi=(\pi \tens \pi)\circ \Delta_C $ and
$\varepsilon_D\circ \pi=\varepsilon_C$.
\smallskip
\par
\textbf{Comodules} : Let \coalg{C} be a coalgebra, a right {\em
$C$-comodule} is an $R$-module $M$ with an $R$-linear map
\rcomod{M} such that the following diagrams are commutative \\
$\xymatrix{
  M \ar[d]_{\roh_M} \ar[r]^{\roh_M} &
  M \tens C \ar[d]^{id\tens\Delta} \\
  M \tens C \ar[r]_(.4){\roh_M\tens id} & M\tens C\tens C}$
\qquad
$ \xymatrix{
  M \ar[dr]_{\cong} \ar[r]^(.4){\roh_M} & M\tens C
                 \ar[d]^{id\tens\varepsilon}  \\
                & M\tens R,             }$ \\
which means, $ (id\tens\Delta)\circ\roh_M=(\roh_M\tens
id)\circ\roh_M $  and  $ (id\tens\varepsilon)\circ\roh_M=id $.
Left $C$-comodules are defined in a symmetric way.
\smallskip
\par
\textbf{Comodule Morphisms} : Let \rcomod{M}, $\roh_N:N\to N\tens
C$ be right $C$-comodules. An $R$-linear map $f:M \to N$ is
called {\em a comodule morphism} (or $C$-{\em colinear}) if
$\roh_N\circ f = (f\tens id)\circ\roh_M $. The set of all
$C$-comodule morphisms between $M$ and $N$ is denoted by
$\comc(M,N)$.
\par
\begin{exam}\label{E.triv}
Let $M$ be a right $C$-comodule, for every $R$-module $W$, $W\tens
M$ has a structure of right $C$-comodule through $W\tens M
\xra{id\tens\roh_M} W\tens M \tens C$. With this structure on
$M\tens C$, $\roh_M$ becomes a right $C$-comodule morphism, which
splits in \rmod. Moreover, for each $R$-linear map $f:W\to V$ the
map $f\tens id:W\tens M \to V\tens M$ is $C$-colinear. The same
can be done analogously for left $C$-comodules.
\end{exam}
\smallskip
\textbf{Bicomodules} : Let $C$, $D$ be two coalgebras, {\em A
$C$-$D$-bicomodule} is a left $C$-comodule and a right
$D$-comodule $M$, such that the left $C$-comodule structure map
$\roh_M^- : M\to C\tens M$ is $D$-colinear, or equivalently the
right $D$-comodule structure map $\roh_M^+ : M\to M\tens D$ is
$C$-colinear.
\smallskip
\par
\textbf{Categories of Comodules} : Let $C$ and $D$ be coalgebras.
The right $C$-comodules together with the $C$-colinear maps
between them constitute an additive category denoted by \comodc.
Analogously, the categories of left $C$-comodules $\ccomod$ and of
$C$-$D$-bicomodules $ {}^C\mathcal M^D $ can be defined.
\begin{lem}\label{L:homcom}
Let $C$ be a coalgebra. Then
\begin{enumerate}[(1)]
 \item The functor $-\tens C:\rmod \to \comodc$ is right adjoint
 to the forgetful functor by the natural isomorphisms (for $M\in
 \comodc, X\in\rmod$),$$\comc(M,X\tens C)\to \homr(M,X),\quad f
\mapsto (id\tens\varepsilon)\circ f.$$
 \item For any $M\in \comodc$, the functor $-\tens M: \rmod \to
 \comodc$ is left adjoint to the functor $\comc(M,-): \comodc \to
 \rmod$ by the natural isomorphisms (for $N\in \comodc, W\in\rmod$),
$$ \comc(W\tens M,N) \to \homr(W,\comc(M,N)),\quad f\mapsto
[w\mapsto f\circ(w\tens -)].$$
\end{enumerate}
\end{lem}
\begin{proof}See~\cite[6.11]{Wis98}
\end{proof}
If $C$ is flat as $R$-module, then the category $\comodc $ is a
Grothendieck category (see~\cite{Wis98}). Moreover a sequence
$$\xymatrix{
  0 \ar[r] & M_1 \ar[r] & M_2 \ar[r] & M_3 \ar[r] & 0 }$$ of
  $C$-comodules is exact in \comodc iff it is exact in \rmod, see
\cite{Sch90}.
\par
\smallskip
{\em From now on we assume that all coalgebras considered in this
work are flat as $R$-modules.}
\section{The Cotensor Functor}\label{s2}
The cotensor functor was first introduced by Milnor and More
\cite{Mil65}  for coalgebras over fields. Guzman \cite{Guzd,Guz89}
studied this functor over rings but in the case of coseparable
coalgebras. One of the key points to recover the Morita-Takeuchi
theorem (see~\cite{Tak77}) in our settings is the associativity of
the cotensor product. In this section we state some results taken
from \cite{Altconf} concerning the exactness and associativity of
the cotensor functor.
\begin{defn}\label{D:defcot}
Let \rcomod{M} be a right and \lcomod{N} a left $C$-comodule. The
cotensor product of $M$ and $N$ (denoted $M\cotensc N$ ) is
defined as the kernel of the $R$-linear map
\[
   \alpha :=\roh_M\tens id_N-id_M\tens \roh_N :\,
    M \tens N \longrightarrow M \tens C \tens N.
\]
So, we have the following exact sequence of $R$-modules
\[
  0\to M \cotensc N \to M \tens N
  \overset{\alpha}{\longrightarrow} M \tens C \tens N
\]
\end{defn}
For a right $C$-comodule $M$, the cotensor functor $M \cotensc- :
\ccomod \to \rmod $ is in general neither left nor right exact
(see \cite{Takh99}), but we have the following

\begin{lem}\label{crexact}
Let $M$ be a right $C$-comodule, then
\begin{enumerate}[(1)]
  \item The functor $M\cotensc - : \ccomod \to \rmod$ is
  $(C,R)$-left exact (also called relative left exact) i.e. left exact
with respect to exact sequences of left $C$-comodules, that are
pure in \rmod.
  \item If $M$ is flat in \rmod, then $M\cotensc -$ is left exact.
\end{enumerate}
\end{lem}
\subsection*{Tensor$-$Cotensor relation}
\begin{lem}\label{L:purtens}
Let $M\in \comodc$, $N\in \ccomod$. For every $W\in\rmod$ there
exist two canonical $R$-linear maps
\begin{eqnarray*}
\gamma_W  : W\tens(M\cotensc N) \lra (W\tens M)\cotensc N, \\
\mu_W  : (M\cotensc N)\tens W \lra M\cotensc (N\tens W).
\end{eqnarray*}
Moreover, the following are equivalent
\begin{enumerate}[(1)]
  \item $M\cotensc N$ is $W$-pure in $M\tens N$.
  \item $\gamma_W$ is an isomorphism.
  \item $\mu_W$ is an isomorphism.
\end{enumerate}
\begin{rem}
$W\tens M$ and $N\tens W$ have the trivial comodule structure
given in example \ref{E.triv}.
\end{rem}
\end{lem}
\subsection*{Associativity of the cotensor product}
One of the main drawbacks of the behaviour of the cotensor product
over rings is that it needs not be associative (see\cite{Gru87}
for a counter example). But under some conditions it becomes so.
Let $M$ be a right $C$-comodule, $L$ a $C$-$D$-bicomodule and $N$
a left $D$-comodule. Then $ M\cotensc L $ has a structure of a
right $D$-comodule through
\[
{M\cotensc L} \xra{id\cotensc\roh^+_L}
 M\cotensc(L\tens D) \cong (M\cotensc L)\tens D,
\]
and it is a subcomodule of $ M\tens L $. Similarly $L\cotensd N $ has a
structure of left $C$-comodule.

Now we give a basic result from \cite{Altconf} that gives necessary
conditions for the associativity of the cotensor product. For
completeness the proof is included.
\begin{prop}[see \cite{Altconf}]\label{P:purcot}
Let $L$ be a $C$-$D$-bicomodule, $M$ a right $C$-comodule and $N$
a left $D$-comodule. If $M\cotensc L$ is $N$-pure in $M\tens L$
and $L\cotensd N$ is $M$-pure in $L\tens N$, then
\[
  M\cotensc(L\cotensd N)\cong (M\cotensc L)\cotensd N.
\]
\end{prop}
\begin{proof}
$M\cotensc L$ and $L\cotensd N$ are right $D$- resp. left $C$-comodules,
since $C$ and $D$ are flat in \rmod. Consider the following commutative
diagram
\[
\xymatrix{
 0 \ar[r] & (M\cotensc L)\cotensd N \ar[r] \ar[d]^{\psi_1} &
(M\tens L) \cotensd N \ar[d]^{\psi_2} \ar[r] & (M\tens C\tens L)\cotensd
N
   \ar[d]^{\psi_3} \\
 0 \ar[r] & M\cotensc(L\cotensd N) \ar[r] & M\tens (L\cotensd N)
\ar[r] & M\tens C\tens (L\cotensd N). }
\]
The rows are exact: The first because $M\cotensc L$ is $N$-pure in
$M\tens L$; the second because it defines $M\cotensc (L\cotensd N)$.
Because of the $M$-purity of $L\cotensd N$ in $L\tens N$ and the
flatness of $C$, $\psi_2$ and $\psi_3$ are isomorphisms. Hence $\psi_1$
is an isomorphism as required.
\end{proof}

\subsection*{Coflat comodules}
\begin{defn}\label{D:coflat}
 A right $C$-comodule $M$ is called {\em coflat} (resp.
{\em faithfully coflat}) if the functor $M\cotensc- :  \ccomod \to
\rmod $ is exact (resp. exact and faithful).
\end{defn}
\begin{prop}\label{eqc.cof}
For a right $C$-comodule $M$ the following are equivalent
\begin{enumerate}
  \item $M$ is faithfully coflat.
  \item The functor $M\cotensc - :\ccomod \lra \rmod$ preserves and
  reflects exact sequences.
  \item $M$ is coflat and $M\cotensc N \neq 0$ for every nonzero left
  $C$-comodule $N$.
\end{enumerate}

\end{prop}
\begin{proof}
This follows from general properties of functors between abelian
categories. (See for example~\cite{Sten})
\end{proof}
If the ground ring is a field, then a $C$-comodule is coflat
(resp. faithfully coflat) iff it is injective (resp. an injective
cogenerator) (see~\cite{Tak277}).
\section{Quasi-finite comodules and the cohom functor}\label{s3}
In this section we investigate the cohom functor, which was first
introduced by M. Takeuchi \cite{Tak77} for coalgebras over fields.
Some technical difficulties arise here, for example we lose the
fact that over fields coflatness and injectivity are equivalent
notions. To overcome this problem some injectivity preserving
condition is imposed. We begin with the definition of
quasi-finite comodules and the cohom functor. Some of the results
given by Takeuchi concerning the cohom functor still hold in our
general settings. In section \ref{s6} we recover all Takeuchi's
results for coalgebras over QF-rings.
\begin{defn} A right $D$-comodule $X$ is called \emph{
quasi-finite}, if the functor $-\tens X : \rmod\to \comodd$ has a
left adjoint. A quasi-finite left $D$-comodule is defined
analogously.
\end{defn}
\begin{nota}
\begin{enumerate}
   \item For a quasi-finite right $D$-comodule X, we denote the left
adjoint of the functor $-\tens X$ by $\cohomf:\comodd \to \rmod
$. This functor is called the cohom functor. We also denote the
adjunction bijection by
\[
  \Phi_{M,W}:\homr(\cohomx{M},W)\to \comd(M,W\tens X).
\]
   \item For every right $D$-comodule $M$ we denote the unit of adjunction
by $\eta_M : M\to \cohomx{M}\tens X$, which satisfies that for
every $D$-colinear map $f:M\to W\tens X$, there exists a unique
$R$-linear map $\tilde{f}:\cohomx{M}\to W$ such that
$f=\Phi_{M,W}(\tilde{f})=(\tilde{f}\tens id)\eta_M$.
\end{enumerate}
\end{nota}
\begin{rems}
\begin{enumerate}
  \item Since the cohom functor has a right adjoint, it follows
that it is right exact and commutes with direct limits.
  \item From the definition it is clear that every quasi-finite
comodule is flat in \rmod.
  \item The image of every $g\in \homr(\cohomx{M},W)$ under the
adjunction bijection can be represented in terms of the counit by
$\Phi_{M,W}(g) = (g\tens id)\eta_M$.
  \item For every right $D$-colinear map $f:N\lra L$, there exists
  a unique $R$-linear map $\cohomx{f} : \cohomx{N}\lra \cohomx{L}$
  such that
\[
\eta_L f=(\cohomx{f}\tens id)\eta_N.
\]
\end{enumerate}
\end{rems}
\subsection*{Some properties of the cohom functor}
\par
Let $X$ be a quasi-finite right $D$-comodule, $M$ a right
$D$-comodule. From the above observations, for every $R$-module
$W$, the $D$-colinear map
\[
id\tens \eta_M : W\tens M \lra W\tens \cohomx{M}\tens X
\]
induces a unique canonical $R$-linear map
\[
\lambda_W : \cohomx{W\tens M} \lra W\tens \cohomx{M}, \]
such that $(\lambda_W\tens id_X)(\eta_{W\tens M}) = id_W\tens
\eta_M$.
\begin{prop}\label{iso.lamb}
The map $\lambda_W$ is an isomorphism.
\end{prop}
\begin{proof}
First it is easy to show that $\lambda_W$ is an isomorphism for
every free $R$-module $W$. Let $W$ be an $R$-module and consider
the free presentation of $W$, $R^{(\Lambda)}\to R^{(\Omega)}\to
W\to 0$. By applying the functors $-\tens M$, $\cohomx{-}$ and
$-\tens \cohomx{M}$ to the free presentation of $W$ we get the
following commutative diagram with exact rows
\[
\xymatrix{ R^{(\Lambda)}\tens \cohomx{M}
\ar@{<-}[d]^{\lambda_{R^{(\Lambda)}}} \ar[r] &
  R^{(\Omega)}\tens \cohomx{M} \ar@{<-}[d]^{\lambda_{R^{(\Omega)}}}
  \ar[r] &
  W\tens  \cohomx{M} \ar@{<-}[d]^{\lambda_W} \ar[r] & 0\\
  \cohomx{R^{(\Lambda)}\tens M} \ar[r] &
  \cohomx{R^{(\Omega)}\tens M} \ar[r] & \cohomx{W\tens M} \ar[r] & 0, }
\]
in which $\lambda_{R^{(\Omega)}}$ and $\lambda_{R^{(\Lambda)}}$
are isomorphisms, hence $\lambda_W$ is an isomorphism.
\end{proof}
\begin{cor}\label{iso.func}
Let $X$ be a quasi-finite right $D$-comodule, $M$ a right
$D$-comodule, then we have the following isomorphism of functors,
\[
\cohomx{-\tens M} \cong -\tens \cohomx{M}:\rmod \to \rmod.
\]
\end{cor}
\begin{proof}
For all $V,W \in \rmod$ we have the natural isomorphisms
\begin{eqnarray*}
   \homr(\cohomx{V\tens M},W)
     &\cong & \comd(V\tens M,W\tens X) \\
     & \cong & \homr(V,\comd(M,W\tens X)) \\
     & \cong & \homr(V,\homr(\cohomx{M},W))\\
     & \cong & \homr(V\tens \cohomx{M},W)
\end{eqnarray*}
\end{proof}
\subsection*{Comodule Structure on $\mathbf{\cohomx{M}}$}
In the following we consider the behaviour of the cohom functor
with respect to bicomodule structures on its arguments. It follows
that the cohom functor is a bifunctor $\cohomb : \bicomod{C}{D}
\times \bicomod{E}{D} \to \bicomod{E}{C}$ contavariant in the
first and  covariant in the second.
\begin{lem}\label{L:dcostr}
\begin{enumerate}[(1)]
\item Let $X$ be a quasi-finite right $D$-comodule, $M$ a
$C$-$D$-bicomodule, then the $R$-linear map
\[
\cohomx{\roh_M^-}:\cohomx{M}\to \cohomx{C\tens M}\cong C\tens
\cohomx{M},
\]
gives $\cohomx{M}$ a structure of a left $C$-comodule. With this
structure the map $\eta_M :M \to \cohomx{M}\tens X$ is
$C$-$D$-bicolinear.
 \item Let $X$ be a $C$-$D$-bicomodule, $M$ a
right $D$-comodule. If $X_D$ is quasi-finite, then there exists an
$R$-linear map
\[
\roh_h: \cohomx{M} \lra \cohomx{M}\tens C
\]
which gives \cohomx{M} a structure of right $C$-comodule,
therefore we have a functor $\cohomf : \comodd \to \comodc$.
Moreover we have $\textrm{Im}(\eta_M)\subseteq \cohomx{M}\cotensc
X$, where $\eta_M$ is the unit of adjunction.
\end{enumerate}
\end{lem}
\begin{proof}
\begin{enumerate}[(1)]
 \item Follows directly from the properties of the unit of adjunction.
 \item The $D$-colinear map $M\xra{(id\tens \roh_X^-)\eta_M}
\cohomx{M}\tens C\tens X$, where $\roh_X^-:X\to C\tens X$ is the
left $C$-comodule structure map, induces the unique $R$-linear map
\[
\roh_h: \cohomx{M}\to \cohomx{M}\tens C,
\]
such that $(id\tens\roh_X^-)\eta_M =(\roh_h\tens id)\eta_M$. This
map gives \cohomx{M} a structure of a right $C$-comodule. To show
that $\cohomf :\comodd \to \comodc$ defines a functor, we have to
show that, for every right $D$-colinear map $f:M\to N$, the
induced map $\cohomx{f} :\cohomx{M} \to \cohomx{N}$ is
$C$-colinear. This follows from the uniqueness of the $R$-linear
map which is induced by the $D$-colinear map $M\xra{(id \tens
\roh_X^-)\circ\eta_{_N}\circ f} \cohomx{N}\tens C\tens X$.

For the other assertion we have
\begin{eqnarray*}
(\roh_h\tens id -id\tens \roh_X^-)\eta_M
   &= & (\roh_h\tens id)\eta_M -(id\tens \roh_X^-)\eta_M \\
   & =& (\roh_h\tens id)\eta_M - (\roh_h\tens id)\eta_M \\
     & = & 0,
\end{eqnarray*}
hence, $\textrm{Im}(\eta_M) \subseteq \textrm{Ker}(\roh_h\tens id
- id\tens \roh_X^-) = \cohomx{M}\cotensc X$.
\end{enumerate}
\end{proof}
The following corollary follows directly from the above
constructions.
\begin{cor}\label{delta}
Let $X$ be a quasi-finite right $D$-comodule.
\begin{enumerate}[(1)]
  \item For every right $D$-comodule $M$ there exists a unique
  $R$-linear map
\[
\delta_M : \cohomx{M} \to M\cotensd \cohomx{D}, \quad
\text{with}\;\; (\delta_M\tens id)\eta_M = id \cotensd \eta_D
\]
  \item If $X$ is also a $C$-$D$-bicomodule, then $\delta_M$ is
right $C$-colinear. If $M$ is a $C$-$D$-bicomodule, then
$\delta_M$ is $C$-bicolinear
\end{enumerate}
\end{cor}
The following theorem was given by Takeuchi \cite{Tak77} for
coalgebras over fields, and it holds also for coalgebras over
rings.
\begin{thm}\label{T:ccohaj}
Let X be a $C$-$D$-bicomodule, then the following are equivalent
\begin{enumerate}
      \item $X_D$ is quasi-finite.
      \item The functor $-\cotensc X:\comodc\to \comodd $ has a
            left adjoint.
      In this case the left adjoint of $-\cotensc X $ is the
         functor $\cohomf : \comodd\to \comodc$.
\end{enumerate}
\end{thm}
\begin{proof}
(1)$\Rightarrow$(2). Define
\begin{alignat*}{1}
\comc(\cohomx{M},N) & \quad \rightleftarrows \quad
\comd(M,N\cotensc X) \\ f & \xymatrix{ \ar@{|->}[r] & }(f\cotensc
id)\eta_M \\
   \tilde{g} & \xymatrix{ \ar@{<-|}[r] & }g,
\end{alignat*}
where $\tilde{g}$ is the unique $R$-linear map
$\tilde{g}:\cohomx{M} \to N$ such that $g=(\tilde{g}\tens
id)\eta_X$. It is straightforward to see that $(f\cotensc
id)\eta_M$ is $D$-colinear, and that the two maps are inverse to
each other. So we have only to show that $\tilde{g}$ is
$C$-colinear. We have
\begin{alignat*}{2}
   (\roh_N \tilde{g}\tens id)\eta_X & = (\roh_N\tens
   id)(\tilde{g}\tens id)\eta_X \\
   & = (\roh_N\tens id)g.
\end{alignat*}
and
\begin{alignat*}{1}
((\tilde{g}\tens id)\roh_h\tens id)\eta_X & = (\tilde{g} \tens
   id\tens id)(\roh_h\tens id)\eta_X \\
   & = (\tilde{g} \tens id\tens id)(id\tens \roh_X^-)\eta_X\\
   & = (\tilde{g}\tens \roh_X^-)\eta_X \\
   & = (id\tens\roh_X^-)(\tilde{g}\tens id)\eta_X \\
   & = (id\tens \roh_X^-)g.
\end{alignat*}
Now $(\roh_N\tens id - id\tens \roh_X^-)|_{N\cotensc X} =0
$, hence $((\tilde{g}\tens id)\roh_h\tens id)\eta_X =
(\roh_N \tilde{g}\tens id)\eta_X $. But an $R$-linear map
$\mu:\cohomx{M}\to N\tens C$, such that $(id\tens \roh_X^-)g =
(\mu\tens id)\eta_X $ is unique, hence $(\tilde{g}\tens
id)\roh_h = \roh_N\tilde{g}$, i.e \~{g} is $C$-colinear.\par
(2)$\Rightarrow$(1). Assume that $-\cotensc X $ has a left
adjoint. The functor $-\tens X : \rmod \to \comodd $ can be
written as the composition of the functors $-\tens C:\rmod\to
\comodc $, which has a left adjoint (see~\ref{L:homcom}), and
$-\cotensc X:\comodc \to \rmod $, this is because for all $W\in
\rmod $, we have $ (W\tens C)\cotensc X \cong  W\tens X $.
\end{proof}
\begin{cor}\label{cunit}
Let $X$ be a $C$-$D$-bicomodule, $X_D$ quasi-finite. For every
$D$-comodule $M$ there exists a $D$-colinear map (the unit of
adjunction) $\eta_M : M\to \cohomx{M}\cotensc X$, such that, for
every $D$-colinear map $f:M\to N\cotensc X$ ($N\in \comodc$),
there exists a unique $C$-colinear map $\tilde{f} : \cohomx{M} \to
N$ with $f=(\tilde{f}\cotensc id)\eta_M$. Moreover, if $M$ is a
$C$-$D$-bicomodule, then $\eta_M$ is $C$-$D$-bicolinear.
\end{cor}
\subsection*{Exactness of the cohom functor}
Now we consider the questions: Under what conditions is the cohom
functor exact; and what happens if this is the case.
\begin{prop}\label{exa.eq}
Let $X$ be a quasi-finite right $D$-comodule. The following are
equivalent
\begin{itemize}
  \item[\emph{(}1\emph{)}] The cohom functor $\cohomf : \comodd \to
  \rmod$ is exact;
  \item[\emph{(}2\emph{)}] $W\tens X$ is injective in $\comodd$, for every
injective $R$-module $W$.
\end{itemize}
If $X$ is a $C$-$D$-bicomodule, then the following are also
equivalent to \emph{(}1\emph{)} and \emph{(}2\emph{)} above.
\begin{itemize}
  \item[\emph{(}3\emph{)}] The cohom functor $\cohomf : \comodd \to
\comodc$ is exact;
  \item[\emph{(}4\emph{)}] $N\cotensc X$ is injective in $\comodd$,
for every injective right $C$-comodule $N$.
\end{itemize}
\end{prop}
\begin{proof}
This is clear, since $\cohomf$ is left adjoint to $-\tens X$ and
the category $\comodd$ has enough injectives.
\end{proof}
Analogous to the definition of injectors in module categories (see
\cite[Exercise 21.7]{And73}), we define injectors in comodule
categories.
\begin{defn}\label{injector}
Let $M$ be a right $C$-comodule. $M$ is said to be an injector in
$\comodc$ if the functor $-\tens M : \rmod \to \comodc$ respects
injective objects. Injectors in $\ccomod$ are defined similarly.
$C$ itself as a comodule is an injector in $\comodc$ and in
$\ccomod$.
\end{defn}
\begin{cor}\label{inj.equ}
Let $X$ be a $C$-$D$-bicomodule with $X_D$ quasi-finite, then the
following are equivalent
\begin{enumerate}[(a)]
  \item $X_D$ is an injector
  \item the cohom functor $\cohomf :\comodd \to \comodc$ is exact
  \item the functor $-\cotensc X :\comodc \to \comodd$ respects
  injective objects.
\end{enumerate}
Now we show that the cohom functor under certain condition is
nothing but a cotensor functor.
\end{cor}
\begin{prop}\label{if.coh.ex}
Let $X$ be a quasi-finite $D$-comodule. If the cohom functor is
exact, then
\begin{enumerate}[(1)]
  \item \cohomx{M} is flat in $\rmod$, for every $D$-comodule $M$,
which is flat in \rmod,
  \item\label{del.is} the $R$-linear map $\delta_M$ in \ref{delta} is an
isomorphism, for every right $C$-comodule $M$,
  \item if, in addition,  $X$ is a $C$-$D$-bicomodule, then the $\delta$
is a bijective $C$-colinear map. If $M$ is also a
$C$-$D$-bicomodule, then $\delta$ is a bijective $C$-bicolinear
map.
\end{enumerate}
\end{prop}
\begin{proof}
\begin{enumerate}[(1)]
  \item We know that $-\tens \cohomx{M} \cong \cohomx{-\tens M}$.
  \item let $M\in \comodd$. The sequence
\[
0 \to M \xra{\roh_M} M\tens D\xra{\roh_M\tens id - id\tens \Delta}
M\tens D\tens D
\]
is exact in \comodd. The result follows now from the following
diagram with exact rows
\[
\xymatrix{
  0 \ar[r] & \cohomx{M} \ar[r]\ar[d]^{\cong} & \cohomx{M\tens D}
  \ar[d]^{\cong}
  \ar[r] & \cohomx{M\tens D\tens D} \ar[d]^{\cong} \\
  0 \ar[r] & M\cotensd \cohomx{D} \ar[r] & M \tens \cohomx{D}
  \ar[r] &M\tens D\tens\cohomx{D}.  }
\]
  \item Follows directly from \ref{del.is}.
\end{enumerate}
\end{proof}
\begin{cor}\label{coh.cot}
Let $X$ be a quasi-finite $D$-comodule. If the cohom functor is
exact then we have
\[
\cohomx{-}\cong -\cotensd \cohomx{D} :\comodd\to \rmod
\]
\end{cor}
\begin{proof}
The isomorphism $\delta_M$ in \ref{if.coh.ex} is natural.
\end{proof}
\section*{The coendomorphism coalgebra}
For a quasi-finite $D$-comodule, and dual to the notion of
endomorphism algebras for modules, the $R$-module $\cohomx{X}$ has
a structure of a coalgebra. This structure was considered in
\cite{Tak77} for fields, and it is also valid in our case.
\vspace{.3cm}
\par
Let $X$ be a quasi-finite $D$-comodule, let $\coend :=\cohomx{X}$.
This $R$-module has a structure of an $R$-coalgebra, where the
comultiplication corresponds to the $D$-colinear map
\[
(id\tens \eta_X)\eta_X : X\to \cohomx{X}\tens \cohomx{X}\tens X
\]
in the adjunction bijection for $M =X$. The counit corresponds to
the $D$-colinear map $id: X\to R\tens X$ for $M=X$ and $W=R$. This
$R$-coalgebra is called the \emph{coendomorphism coalgebra} of
$X$. The unit of adjunction $\eta_X:X\to \cohomx{X}\tens X$ gives
$X$ a structure of a left $\coend$-comodule, hence $X$ is a
$\coend$-$D$-coalgebra.

\vspace{.3cm}
\par
It is well known that for an $R$-coalgebra \coalg{C}, the dual
$R$-module $C^\ast =\homr(C,R)$ is an $R$-algebra under the
convolution product( i.e for $f,g \in C^\ast,\, f\ast g =(f\tens
g)\Delta )$. We now give a result about the dual algebra of the
coendomorphism coalgebra.
\begin{prop}\label{Pr:anti}
For a quasi-finite right $D$-comodule $X$, the dual algebra of the
coendomorphism coalgebra is algebra anti-isomorphic to the algebra
of right $D$-comodule endomorphisms of $X$. (i.e.
$(\coend)^{\ast} \cong \comd(X,X)$).
\end{prop}
\begin{proof}
Let $X$ be a quasi-finite right $D$-comodule, and consider the
composition of the maps $$ \coend^\ast =
\homr(\cohomx{X},R)\xra[\cong]{\Phi_{X,R}} \comd(X,R\tens X)\cong
\comd(X,X).$$ The image of $f\in \coend^\ast$ under this
composition is $\vartheta (f\tens id)\circ\eta_X$, where
$\vartheta:R\tens X \to X$ is the canonical isomorphism. It is
easy to see that this map is an algebra anti-isomorphism.
\end{proof}
\begin{rem}
In \ref{Pr:anti}, if we start with \emph{a left} quasi-finite
$D$-comodule, then the algebra anti-isomorphism there is an
algebra isomorphism.
\end{rem}
\section{Functors between comodule categories}\label{s4}
Let $C, D$ be $R$-coalgebras, we will study additive functors
$F:\comodc \to \comodd $, which are assumed to be $R$-linear, i.e
the canonical maps $ \comc(M,N) \to \comd(F(M),F(N)) $, induced by
$F$ for $M, N \in \comodc $, are assumed to be $R$-linear. We will
show that good enough $R$-linear functors are isomorphic to a
cotensor functor.
Let $F:\comodc \to \comodd $ be an $R$-linear functor. We consider
the bifunctors
\[
-\tens F(-),\; F(-\tens -):\rmod \times \comodc \to \comodd,
\]
where for $W\in \rmod,\; M\in \comodc $, $W\tens M$ (resp. $W\tens
F(M)$) is endowed with the canonical structure of right
$C$-comodule (resp. $D$-comodule). Our aim is to construct a
natural transformation
\[
\lambda : -\tens F(-) \to F(-\tens -).
\]
For every $C$-comodule $M$ let $\lamb{R}{M}$ be the unique
isomorphism, that render the following diagram commutative
\begin{equation*}
 \xymatrix{
  R\tens F(M) \ar[d]_{\cong} \ar[r]^{\lamb{R}{M}}
                & F(R\tens M) \ar[d]^{\cong}  \\
  F(M) \ar[r]_{=}
                & F(M),             }
\end{equation*}
where $M \in \comodc $ and the vertical arrows are the canonical
isomorphisms.
\begin{lem}
Let $\mathcal{P}=\{R\} $ be the full subcategory of $\rmod$ whose
only object is $R$, denote the restriction of the functors $-\tens
F(M)$ and $F(-\tens M)$ on $\mathcal{P}$ by $T$ and $S$ resp. Then
$\lamb{-}{M} : T\to S$ is a natural transformation.
\end{lem}
\begin{proof}
Given a homomorphism $f:R\to R$ in $\mathcal{P}$, we have to
check that the diagram $$ \xymatrix{
  R\tens F(M) \ar[d]_{\lamb{R}{M}} \ar[r]^{f\tens id} & R\tens F(M)
\ar[d]^{\lamb{R}{M}} \\ F(R\tens M) \ar[r]_{F(f\tens id)} &
F(R\tens M)} $$ commutes.

Define $g:F(M)\to F(M)$, by $g(x)=f(1)x$. It is easy to see that
$g$ is $D$-colinear and that $g=F(\tilde{g})$, where $\tilde{g} :
M \to M, m\mapsto f(1)m$.

Consider the diagram
\begin{equation*}
\xymatrix{
   & F(M) \ar[rr]^{g} \ar'[d][dd]_{id}
   & & F(M) \ar[dd]^{id}        \\
  R\tens F(M) \ar[ur]^{\cong}\ar[rr]^(.7){f\tens id}\ar[dd]_
  {\lamb{R}{M}} & & R\tens F(M) \ar[ur]_{\cong}\ar[dd]^(.3){\lamb{R
  }{M}}  & \\
  & F(M) \ar'[r]^(.6){g}[rr]
      & & F(M)                \\
  F(R\tens M) \ar[rr]_{F(f\tens id)} \ar[ur]^{\cong}
      & & F(R\tens M) \ar[ur]_{\cong}        }
\end{equation*}
By the definition of $\lamb{R}{M}$ the right and left rectangles
are commutative, the upper rectangle is commutative because
$f(r)=f(1)r$ for all $r\in R$. The lower rectangle is
commutative, since it is obtained from the diagram
\[
\xymatrix{
  M \ar[r]^{\tilde{g}}
                & M \\
  R\tens M  \ar[r]_{f\tens id} \ar[u]^{\cong}
                & R\tens M \ar[u]^{\cong}    }
\]
which is commutative. Now, since the map $F(R\tens M) \to F(M)$
is a monomorph (isomorph), the front rectangle is commutative.
This implies that $\lamb{-}{M}:-\tens F(M) \to F(-\tens M) $ is
natural, where $-\tens F(M),\; F(-\tens M):\mathcal{P}\to
\comodd$.
\end{proof}
By a theorem of Mitchel (\cite[Theorem 3.6.5]{Pop73}), \lamb{-}{M}
is uniquely extended to a natural transformation
\begin{equation*}
\lam{-}{M}:-\tens F(M) \to F(-\tens M),
\end{equation*}
for every right $C$-comodule $M$, where now $-\tens F(M), \;
F(-\tens M) : \rmod \to \comodd$. Moreover (see \cite[Corollary
3.6.6]{Pop73}), if $F$ preserves direct sums (resp. direct limits,
inductive limits), then \lam{W}{-} is a natural isomorphism for
every projective $R$-module (resp. flat $R$-module, $R$-module)
$W$.
Now we want to see under what conditions is $\lam{W}{M}$
functorial in $M$.
\begin{lem}
Assume that $F:\comodc \to \comodd$ respects direct sums, then
\[
\lam{-}{-} :-\tens F(-) \to F(-\tens -)
\]
is a natural transformation. Moreover, if $F$ preserves direct
limits, then $\lam{W}{-}$ is a natural isomorphism for every flat
$R$-module $W$. Finally, if $F$ preserves inductive limits, then
$\lam{-}{-}$ is a natural isomorphism.
\end{lem}
\begin{proof}
We know that $\lam{W}{M}$ is functorial in $W$. Thus we have to
show that it is functorial in $M$. Let $g:M\to N$ be $C$-colinear.
From the commutative diagram
\[
   \xymatrix{
   R\tens F(M) \ar[rr]^{id\tens F(f)} \ar[dd]_{\lam{R}{M}}
    \ar[dr]^{\cong}  & &  R\tens F(N) \ar[dd]^(.6){\lam{R}{N}}
    \ar[rd]^{\cong}  &      \\
  & F(M) \ar'[r]^(.6){F(f)}[rr]
      & &  F(N)  \\
   F(R\tens M) \ar[rr]_{F(id\tens f)} \ar[ru]^{\cong}
      & & F(R\tens N) \ar[ur]_{\cong} &       }
\]
we obtain that $\lam{R}{M}$ is functorial in $M$. Now, since $F$
preserves direct sums, $\lam{W}{M}$ is natural for every free
$R$-module $W$. In the general case use a free presentation for
$W$ to show that \lam{W}{M} is natural (in $M$) for every
$R$-module $W$. The other assertions follow from the above
observations.
\end{proof}
\section*{Bicomodule structure on $\mathbf{F(C)}$}
Now, using the natural transformation $\lambda$ and by imposing a
mild condition on $F$, we want to give $F(C)$ a structure of a
bicomodule. To this end we need the following lemma
\begin{lem}\label{formula}
Assume that $F:\comodc\to\comodd$ preserves direct sums, then for
every $W\in \rmod$ and $M, N \in \comodc$, the following formula
holds
\[
\lam{W}{M\tens N}\circ (id\tens \lam{M}{N}) = \lam{W\tens M}{N}
\]
\end{lem}
\begin{proof}
It is easy to prove the result for $W=R$ and hence for every free
$R$-module $W$. Now let $W \in \rmod$, from the free presentation
$R^{(\Lambda)}\xra{\;\;\theta\;\;} W\to 0$ we obtain the following
diagram
\[
{\mbox {\footnotesize \xymatrix{ R^{(\Lambda)}\tens F(M\tens N)
\ar[dd]_{\lam{R^{(\Lambda)}}{M\tens N}} \ar[rr]^{\theta\tens
F(id)} \ar@{<-}[dr]^{id\tens \lam{M}{N}}&  & W\tens F(M\tens N)
\ar[dd]_(.8){\lam{W}{M\tens N}} \ar@{<-}[dr]^{id\tens \lam{M}{N}}
&
\\ & R^{(\Lambda)}\tens M \tens F(N)
\ar@{->>}'[r]^(.7){\theta\tens id}[rr] & & W\tens M\tens F(N) \\
F(R^{(\Lambda)}\tens M\tens N)
\ar@{<-}[ur]^{\lam{R^{(\Lambda)}\tens M}{N}}
\ar[rr]_{F(\theta\tens id)} & & F(W\tens M\tens N),
\ar@{<-}[ur]_{\lam{W\tens M}{N}}}}}
\]
in which all sub-diagrams (except possibly the right triangle) are
commutative and $\theta \tens id$ is an epimorph. Hence the right
triangle is commutative and the result follows.
\end{proof}
\begin{prop}\label{bicof(c)}
Assume that $F$ preserves direct limits, then $F(C)$ has a
structure of a left $C$-comodule. Hence $F(C)$ becomes a
$C$-$D$-bicomodule.
\end{prop}
\begin{proof}
Since $F$ preserves direct limits, $\lam{C}{C}$ is an isomorphism,
since $C_R$ is flat. Let $ \roh_{F(C)} : F(C)\lra C\tens F(C)$ be
the unique $R$-linear map that makes the diagram
\[
\xymatrix{
   F(C) \ar[rr]^{\roh_{F(C)}} \ar[dr]_{F(\Delta)} & & C\tens
   F(C)\ar[dl]_{\cong}^{\lam{C}{C}} \\
   & F(C\tens C)     }
\]
commutative.

We will show that $\roh_{F(C)}$ is a comodule structure map.
First, we have to show that the following diagram is commutative
\[
 \xymatrix{
  F(C) \ar[d]_{\roh_{F(C)}} \ar[r]^{\roh_{F(C)}}
  \ar@{}[rd]^{(\star)}
  & C\tens F(C) \ar[d]^{\Delta\tens id}  \\
  C\tens F(C)  \ar[r]_{id\tens \roh_{F(C)}}
                & C\tens C\tens F(C).             }
\]
Imbed this diagram in the following diagram
\[
\xymatrix{
   & F(C) \ar@{->}'[d]'[ddd]^{id}[dddd] \ar[rr]^{\roh_{F(C)}}
   \ar[dl]_{\roh_{F(C)}} &&  C\tens F(C)
   \ar[dddd]^{\lam{C}{C}} \ar[dl]_{\Delta\tens id}\\
  C\tens F(C) \ar[dddd]_{\lam{C}{C}}
    \ar[rdd]|{id\tens F(\Delta)} \ar[rr]^(.64){id\tens \roh_{F(C)}} &&
    C\tens C\tens F(C) \ar[dddd]|(.4){\lam{C\tens C}{C}} \ar[ddl]|
    {id\tens \lam{C}{C}} &  \\ &&&\\
    &  C\tens F(C\tens C) \ar[rdd]|(.65){\lam{C}
 {C\tens C}} &&  \\  & F(C) \ar[rr]^(.45){F(\Delta)} \ar[dl]^{F(\Delta)} &
  & F(C\tens C) \ar[dl]^{F(\Delta\tens id)} \\
  F(C\tens C)  \ar[rr]_{F(id\tens\Delta)} & &
    F(C\tens C\tens C)  &    }
\]
This diagram (except possibly the top rectangle) is commutative,
this is shown using lemma \ref{formula}, the coassociativity of
$\Delta$, that \lam{C\tens C}{C} is an isomorphism, and that
\lam{-}{C} is natural. Hence the top side is commutative since
\lam{C\tens C}{C} is a monomorphism.

To show that $(\varepsilon \tens id)\roh_{F(C)} = id$, consider
the diagram $$ \xymatrix{
  R\tens F(C) \ar[rr]^{\lam{R}{C}} &&
  F(R\tens C)  \\
  & F(C)\ar[dr]^{F(\Delta)} \ar[ru]^{\cong} \ar[ul]_{\cong}
  \ar[ld]_{\roh_{F(C)}} & \\
  C\tens F(C) \ar[uu]^{\varepsilon\tens id} \ar[rr]_{\lam{C}{C}} &&
  F(C\tens C) \ar[uu]_{F(\varepsilon\tens id)}   }$$
which is shown to be commutative by the counitary property, and
the definition of $\roh_{F(C)}$.
\end{proof}
We know that the cotensor functor is relative left exact and
respects direct limits. Now we give a theorem, similar to Watts
theorem for modules, which shows that a functor that is relative
left exact and respects direct limits with an extra condition can
be presented as a cotensor functor.
\begin{thm}
Let $F:\comodc\to \comodd$ be a relative left exact functor that
respects direct limits. If $\lam{M}{C}$ and $\lam{M\tens C}{C}$
are isomorphisms, for every right $C$-comodule $M$ (e.g. if $F$
respects inductive limits or $F$ is an equivalence), then $F$ is
naturally isomorphic to $-\cotensc F(C)$.
\end{thm}
\begin{proof}
Let $\rcomod{M}$ be a right $C$-comodule. We have the following
$(C,R)$-exact sequence of comodules
\[
0\lra M\lra M\tens C\lra M\tens C\tens C,
\]
hence we get the following commutative (see \ref{formula}) diagram
with exact rows
\[
\xymatrix{
  0 \ar[r] & M\cotensc F(C) \ar[r] \ar@{-->}[d]^{\phi_M} &
  M\tens F(C) \ar[rrr]^(.45){\roh_M\tens id-id\tens \roh_{F(C)}}
  \ar[d]^{\lam{M}{C}} & & &
  M\tens C\tens F(C) \ar[d]^{\lam{M\tens C}{C}} \\
 0\ar[r]& F(M) \ar[r] & F(M\tens C) \ar[rrr]_(.45){F(\roh_M\tens id-id
\tens\Delta)} & & & F(M\tens C\tens C),  }
\]
where the desired isomorphism $\phi_M$ is given by the universal
property of the kernel. To show the naturality of $\phi_M$, let
$f:M\to N$ be a $C$-colinear map and use the following commutative
diagram
\[
\xymatrix{
  &  M\tens F(C) \ar[rr]^{f\tens id} \ar'[d][dd]_{\lam{M}{C}}
       & & N\tens F(C) \ar[dd]^{\lam{N}{C}}        \\
       M\cotensc F(C) \ar[rr]^(.7){f\cotensc id} \ar@{>->}[ur]
       \ar[dd]_(.7){\phi_M}
    & & N\cotensc F(C) \ar@{>->}[ur] \ar[dd]_(.7){\phi_N}& \\
    & F(M\tens C)\ar'[r][rr]^(.4){F(f\tens id)}
       & &  F(N\tens C) \\
F(M)\ar@{>->}[ur]^{F(\roh_M)}\ar[rr]_{F(f)}
  & & F(N) \ar@{>->}[ur]_{F(\roh_N)} & }
\]
to see that $\phi_M$ is functorial in $M$.
\end{proof}
\begin{cor}\label{c.equiv}
Let $C$, $D$ be $R$-coalgebras. If $F :\comodc \to \comodd$ is a
category equivalence with inverse $G: \comodd\to\comodc$, then
there exist bicomodules $X\in \cdbicomod$ and $Y\in \dcbicomod$
such that
\[
F\cong -\cotensc X \quad \text{and} \quad G\cong -\cotensd Y.
\]
\end{cor}
\begin{proof}
$X=F(C)$ and $Y=G(D)$.
\end{proof}
\section{Equivalences of Comodule Categories}\label{s5}
In this section we study equivalences between comodule categories,
and under what conditions we get such an equivalence. We prove a
generalization of Morita-Takeuchi theorem for our settings. We
begin with some properties of the bicomodules $F(C)$ and $G(D)$ in
corollary \ref{c.equiv}.
\begin{lem}\label{inv.equ}
Let $F :\comodc \to \comodd$ be a category equivalence with
inverse $G: \comodd\to\comodc$. Then the following hold
\begin{enumerate}[(1)]
  \item $F(C)$ and $G(D)$ are flat $R$-modules.
  \item For each right $C$-comodule $M$, $M\cotensc F(C)$ is a
  pure submodule of $M\tens F(C)$.
  \item For each right $D$-comodule $N$, $N\cotensd G(D)$ is pure in
  $N\tens G(D)$.
  \item For every $P\in \text{\dcomod}$ we have
\[
  \mathbf{\bigl(} G(D)\cotensc F(C)\mathbf{\bigr)}\cotensd P \cong
G(D)\cotensc \mathbf{\bigl(}F(C)\cotensd P\mathbf{\bigr)}
\]
  \item For every $Q \in \text{\ccomod}$ we have
\[
\mathbf{\bigl(}F(C)\cotensd G(D)\mathbf{\bigr)}\cotensc Q\cong
F(C)\cotensd \mathbf{\bigl(}G(D)\cotensc Q\mathbf{\bigr)}.
\]
\end{enumerate}
\end{lem}
\begin{proof}
\begin{enumerate}[(1)]
    \item $-\tens F(C) \cong F(-\tens C)$ which is exact, since
    $C$ is flat in \rmod and $F$ is an equivalence. Similarly for
    $G(D)$.
    \item\label{r} Let $M\in \comodc$. For every  $W\in \rmod$ we
    have the canonical $R$-linear map $\gamma_W : W\tens (M\cotensc F(C))
\to (W\tens M)\cotensc F(C)$, (see lemma \ref{L:purtens}).
Consider the following commutative diagram
\[
\xymatrix{
     0 \ar[r] & F(W\tens M) \ar[r] \ar[d]^{\cong} & F(W\tens M\tens C)
  \ar[r]\ar[d]^{\cong} & F(W\tens M\tens C\tens C) \ar[d]^{\cong}  \\
     0\ar[r] & W\tens (M\cotensc F(C)) \ar[r] \ar[d]^{\gamma_W} &
  W\tens (M\tens F(C)) \ar[d]^{\cong} \ar[r] & W\tens (M\tens
  C\tens F(C)) \ar[d]^{\cong} \\
     0\ar[r] &(W\tens M)\cotensc F(C) \ar[r] & (W\tens M)\tens
     F(C)\ar[r] & (W\tens M)\tens C\tens F(C), }
\]
in which the first and third rows are exact, hence the second is
also exact. The result follows now from lemma \ref{L:purtens}.

   \item Similar to \ref{r}.
    \item The assertions in (4) and (5) follow from (1), (2)
   and proposition \ref{P:purcot}.
\end{enumerate}
\end{proof}

With the help of representing an equivalence by a cotensor
functor and lemma \ref{inv.equ}, we can now prove that an
equivalence between right comodule categories gives an
equivalence of the corresponding left comodule categories
\begin{thm}
Let $C, D$ be two coalgebras. If $\comodc$ is equivalent to
$\comodd$, then $\ccomod$ is equivalent to $\dcomod$.
\end{thm}
\begin{proof}
Let $F:\comodc \to \comodd$ be an equivalence with inverse
$G:\comodd\to \comodc$. Define $F':=F(C)\cotensd -:\text{\dcomod}
\to \text{\ccomod}$ and $G':=G(D)\cotensc -:\text{\ccomod} \to
\text{\dcomod}$. For $P\in$ \dcomod, $Q\in$ \ccomod, we have
\begin{eqnarray*}
F'G'(Q) & = & F(C)\cotensd \mathbf{\bigl(}G(D) \cotensc Q\mathbf{\bigr)} \\
        & \cong & \mathbf{\bigl(}F(C) \cotensd G(D)\mathbf{\bigr)}
                                 \cotensc Q \\
        & \cong & GF(C)\cotensc Q \\
        & \cong & C\cotensc Q \cong Q. \\
G'F'(P) & = & G(D)\cotensc \mathbf{\bigl(}F(C) \cotensd P\mathbf{\bigr)} \\
        & \cong & \mathbf{\bigl(}G(D)\cotensc F(C)\mathbf{\bigr)}
                    \cotensd P \\
        & \cong & FG(D) \cotensd P \\
        & \cong & D\cotensd P \cong P.
\end{eqnarray*}
Hence $F'$ is an equivalence with inverse $G'$
\end{proof}
\begin{defn}
Two coalgebras are called Morita-Takeuchi equivalent if the
categories of right (equivalently of left) comodules over these
coalgebras are equivalent
\end{defn}
The bicomodules $F(C)$ and $G(D)$ have other properties
\begin{cor}\label{c.qf.coflat}
Let $F:\comodc\to \comodd$ be a category equivalence with inverse
$G:\comodd\to \comodc$. Then $F(C)$ is  quasi-finite, faithfully
coflat and an injector as right $D$- and as left $C$-comodule.
Moreover $e_D(F(C)) \cong C$ and $e_C(F(C)) \cong D$ as
coalgebras. Similar results hold for $G(D)$.
\end{cor}
\begin{proof}
The fact that $F\cong -\cotensc F(C)$ is an equivalence with
inverse $G\cong -\cotensd G(D)$ implies that $F(C)$ is
quasi-finite and an injector as right $D$-comodule and faithfully
coflat as left $C$-comodule, and that $e_D(F(C)) \cong
F(C)\cotensd G(D) \cong GF(D) \cong C$. The other assertions
follow from the fact that $F(C)\cotensd -$ is an equivalence with
inverse $G(D)\cotensc -$.
\end{proof}
We are now ready to give our main result in this article, which
generalizes the Morita-Takeuchi theorem.
\begin{thm}\label{mtt}
Let $C$ and $D$ be two coalgebras. The following are equivalent.
\begin{enumerate}[(1)]
  \item\label{e1} $C$ and $D$ are Morita-Takeuchi equivalent.
  \item The categories $\ccomod$ and $\dcomod$ are equivalent.
  \item\label{e2} There exists a $C$-$D$-bicomodule $X$, such that $X_D$ is
  quasi-finite, faithfully coflat and an injector, and $e_D(X)\cong
  C$ as coalgebras.
  \item There exists a $C$-$D$-bicomodule $X$, such that ${}_CX$ is
  quasi-finite, faithfully coflat and an injector, and $e_C(X)\cong
  D$ as coalgebras.
  \item There exists a $D$-$C$-bicomodule $Y$, such that $Y_C$ is
  quasi-finite, faithfully coflat and an injector, and $e_C(X)\cong
  D$ as coalgebras.
  \item There exists a $D$-$C$-bicomodule $Y$, such that ${}_DY$ is
  quasi-finite, faithfully coflat and an injector, and $e_D(Y)\cong
  C$ as coalgebras.
\end{enumerate}
\end{thm}
\begin{proof}
From what we have done, we have  only to show that (\ref{e2})
$\Rightarrow$ (\ref{e1}). Assume that a $C$-$D$-bicomodule $X$
satisfies the conditions of (\ref{e2}). We show that $-\cotensc
X:\comodc\to \comodd$ is an equivalence.

From corollary \ref{coh.cot} we know that $\cohomx{-}\cong
-\cotensd \cohomx{D}$, hence $\cohomx{D}$ is a coflat
$D$-comodule. Next we show that $\eta_D:D\to\cohomx{D}\cotensc X$,
(see \ref{cunit}), is an isomorphism. The $C$-bicolinear map
$\delta_X:\cohomx{X}\cong C \to X\cotensd \cohomx{D}$ is an
isomorphism (see \ref{if.coh.ex}) and the diagram
\[
\xymatrix{
  X \ar[d]_{\cong} \ar[rr]^{\cong} &
                & X\cotensd D \ar[d]^{id_X\cotensd \eta_D}  \\
  C\cotensc X  \ar[rr]_(.35){\delta_X\cotensc id_X}
              &  & X\cotensd\cohomx{D}\cotensc X }
\]
is commutative (Notice that $X\cotensd \Bigl(\cohomx{D} \cotensc
X\Bigr)\cong \Bigl(X\cotensd \cohomx{D}\Bigr)\cotensc X$, since
$X_D$ is coflat). Hence $id_X\cotensd \eta_D$ is an isomorphism,
therefore, since $X_D$ is faithfully coflat, $\eta_D$ is an
isomorphism.

Finally we show that $\cohomx{D}\cotensc X$ is pure in
$\cohomx{D}\tens X$. Let $W\in \rmod$ and consider the following
commutative diagram, where $\mu_W$ is the canonical map of lemma
\ref{L:purtens}%
\[
{\mbox {\scriptsize \xymatrix{
  0  \ar[r] & X\cotensd \Bigl(\cohomx{D}\cotensc (X\tens W)\Bigr)
   \ar[r] \ar@{<-}[d]^{id_X \cotensc \mu_W}
   & X\cotensd \Bigl(\cohomx{D}\tens X \tens W\Bigr)
    \ar@{=}[d]
  \ar[r] & X\cotensd \Bigl(\cohomx{D}\tens C\tens X \tens W\Bigr)
    \ar@{=}[d] \\
   0\ar[r] & X\cotensd \Bigl((\cohomx{D}\cotensc X)\tens W\Bigr)
   \ar[d]^{\cong} \ar[r] & X\cotensd \Bigl(\cohomx{D}\tens X
   \tens W\Bigr) \ar[d]^{\cong}
  \ar[r] & X\cotensd \Bigl(\cohomx{D}\tens C\tens X \tens W\Bigr)
  \ar[d]^{\cong} \\
  0 \ar[r] & \Bigl(X\cotensd \cohomx{D}\Bigr) \cotensc (X\tens W)
   \ar[r] & \Bigl(X\cotensd \cohomx{D}\Bigr)\tens (X\tens W)
  \ar[r] & \Bigl(X\cotensd \cohomx{D}\Bigr)\tens C\tens X\tens
  W,}}}
\]
in which the first and third rows are exact, since $X_D$ is
coflat. So the second row is also exact and hence $id_X\cotensd
\mu_W$ is an isomorphism. From lemma (\ref{L:purtens}) it follows
that $\cohomx{D}\cotensc X$ is pure in $\cohomx{D}\tens X$.

Now let $M\in \comodd$ and $N\in \comodc$, we have
\begin{eqnarray*}
\cohomx{N\cotensc X} & \cong & (N\cotensc X)\cotensd \cohomx{D} \\
           & \cong & N\cotensc (X\cotensd \cohomx{D}) \quad
             \text{(since $\cohomx{D}$ is coflat in \ccomod)}\\
           &\cong & N\cotensc \cohomx{X} \\
           &\cong & N\cotensc C \cong N, \quad \text{and}\\
\cohomx{M}\cotensc X &\cong & (M\cotensd \cohomx{D})\cotensc X \\
           & \cong & M\cotensd (\cohomx{D}\cotensc X)\\
           &\cong & M\cotensd D \cong M.
\end{eqnarray*}
Therefore $-\cotensc X$ is an equivalence with inverse $-\cotensd
\cohomx{-}$.
\end{proof}
\begin{defn}\label{invertible}
A $C$-$D$-bicomodule $M$ is called invertible if the functor
$-\cotensc M :\comodc \to\comodd$ is an equivalence.
\end{defn}

\begin{cor}\label{inv.eq}
Let $C$, $D$ be two $R$-coalgebras. For a $C$-$D$-bicomodule $X$,
the following are equivalent
\begin{enumerate}[(1)]
  \item $X$ is invertible.
  \item $X\cotensd -: \dcomod \to \ccomod$ is an equivalence.
  \item $X_D$ is quasi-finite, faithfully coflat and an injector,
  and $e_D(X) \cong C$ as coalgebras.
  \item ${}_CX$ is quasi-finite, faithfully coflat and an injector,
  and $e_C(X) \cong D$ as coalgebras.
\end{enumerate}
In this case, $C$ and $D$ are Morita-Takeuchi equivalent.
\end{cor}
\section{Morita-Takeuchi Context}\label{Scontext}
In this section we give a definition of Morita-Takeuchi context
for coalgebras over rings. We also establish a correspondence
between equivalences of comodule categories and strict
Morita-Takeuchi contexts.
\begin{defn}
A \emph{Morita-Takeuchi context} ($D, C, M, N, f,g$) consists of
$R$-coalgebras $D$ and $C$, bicomodules ${}_DM_C$, ${}_CN_D$, and
bicomodule morphisms $f:D\lra M\cotensc N$, $g:C\lra N\cotensd M$,
such that the following conditions hold
\begin{enumerate}[(1)]
  \item\label{co1} $M$ and $N$ are flat as $R$-modules.
  \item\label{co2} $M\cotensc N$ resp. $N\cotensd M$ are pure in $M\tens N$
  resp. $N\tens M$.
  \item The diagrams
\[
\xymatrix{
  M \ar[d]_{\cong} \ar[r]^{\cong} & M\cotensc C \ar[d]^{id\cotensc g}  \\
  D\cotensd M  \ar[r]_(.4){f\cotensd id} & M\cotensc N\cotensd M,} \quad
\xymatrix{
  N \ar[d]_{\cong} \ar[r]^{\cong} & N\cotensd D \ar[d]^{id\cotensd f}  \\
  C\cotensc N  \ar[r]_(.4){g\cotensc id} & N\cotensc M\cotensd N.}
\]
commute.
\end{enumerate}
\end{defn}
\begin{rems}
\begin{enumerate}[(a)]
  \item Caenepeel \cite{Cae98} defined a Morita-Takeuchi context
  for coalgebras over rings which is the same as our definition
  but without the purity condition, with his definition it is not
  possible to prove that a strict Morita-Takeuchi context gives an
  equivalence between the involved coalgebras, since in his case
  the cotensor product is not associative.
  \item If $R$ is a field, then conditions \ref{co1} and \ref{co2}
are satisfied for all $M$ and $N$. In this case our definition
reduces to Takeuchi's one (see \cite{Tak77}).
  \item From condition \ref{co2} it follows that
\begin{eqnarray*}
 X\cotensc (N\cotensd M) &\cong & (X\cotensc N)\cotensd M, \qquad
   Y\cotensd (M\cotensc N)\;\; \cong \;\; (Y\cotensd M)\cotensc N,   \\
 N\cotensd (M\cotensc Z) & \cong & (N\cotensd M)\cotensc Z,  \qquad
   M\cotensc (N\cotensd T) \;\; \cong \;\; (M\cotensc N)\cotensd T,
\end{eqnarray*}
for all $X\in \comodc$, $Y\in \comodd$, $Z\in\ccomod$ and
$T\in\dcomod$.
\end{enumerate}
\end{rems}
From a Morita-Takeuchi context ($D, C, M, N, f,g$) we have the
following left exact functors
\begin{eqnarray*}
 F &:=-\cotensd M :\comodd \lra \comodc, \\
 G &:=-\cotensc N : \comodc \lra \comodd, \\
 F' &:= M\cotensc -: \text{\ccomod} \lra \dcomod,  \\
 G' &:= N\cotensd -  :\text{\dcomod} \lra \ccomod,
\end{eqnarray*}
and the following natural transformations
\begin{eqnarray*}
\lambda : id_{\comodc} \to FG, \quad & \zeta : id_{\comodd} \to
GF, \\ \psi : id_{\text{\dcomod}}\to F'G' \quad  & \delta :
id_{\text{\ccomod}}\to G'F',
\end{eqnarray*}
which are defined as follows
\begin{eqnarray*}
 \lambda_X &: X\cong X\cotensc C \xra{id\cotensc g} X\cotensc
    (N\cotensd M) \cong (X\cotensc N)\cotensd M = FG(X),\\
 \zeta_Y &: Y \cong Y\cotensd D \xra{id\cotensd f} Y\cotensd
    (M\cotensc N)\cong (Y\cotensd M)\cotensc N = GF(Y), \\
 \psi_T &: T \cong D\cotensd T \xra{f\cotensd id} (M\cotensc N)
    \cotensd T \cong M\cotensc (N\cotensd T) = F'G'(T),\\
 \delta_Z &: Z\cong C\cotensc Z \xra{g\cotensc id} (N\cotensd M)
    \cotensc Z \cong N\cotensd (M\cotensc Z) = G'F'(Z),
\end{eqnarray*}
\begin{lem}\label{f.inj}
Let ($D, C, M, N, f,g$) be a Morita-Takeuchi context.
\begin{enumerate}[(1)]
  \item If $f$ is bijective, then
\begin{enumerate}[(a)]
  \item $G=-\cotensc N : \comodc \lra \comodd$ is left
  adjoint to $F =-\cotensd M :\comodd \lra \comodc$.
  \item $F'=M\cotensc - : \ccomod \lra \dcomod$ is left
  adjoint to $G' =N\cotensd - :\dcomod \lra \ccomod$.
\end{enumerate}
  \item If $g$ is bijective, then
\begin{enumerate}[(a)]
  \item $F =-\cotensd M :\comodd \lra \comodc$ is left
  adjoint to $G=-\cotensc N : \comodc \lra \comodd$
  \item $G'=N\cotensd - : \dcomod \lra \ccomod$ is
  left adjoint to $F'=M\cotensc - : \ccomod \lra \dcomod$.
\end{enumerate}
  \item If $f$ and $g$ are bijective, then
\begin{enumerate}[(a)]
  \item $G=-\cotensc N$ is an equivalence
  with inverse $F =-\cotensd M$.
  \item $G'=N\cotensd -$ is an equivalence
  with inverse $F'=M\cotensc -$.
\end{enumerate}
In this case $C$ and $D$ are Morita-Takeuchi equivalent
coalgebras.
\end{enumerate}
\end{lem}
\begin{proof}
\begin{enumerate}[(1)]
 \item Assume $f$ is bijective
\begin{enumerate}[(a)]
  \item\label{F.ladj} We have the natural transformation
$\lambda : id_{\comodc} \lra FG$. Since $f$ is bijective
$\zeta^{-1} : GF \lra id_{\comodd}$ is a natural transformation.
We are done if we show that each of the composition of morphisms
\begin{eqnarray*}
& & X\cotensc N \cong X\cotensc C\cotensc N \xra{id\cotensc
g\cotensc id} X\cotensc N\cotensd M\cotensc N \xra{id\cotensc
id\cotensd
f^{-1}} X\cotensc N, \\
& & Y\cotensd M \cong Y\cotensd M\cotensc C \xra{id\cotensd
id\cotensc g} Y\cotensd M\cotensc N\cotensd M \xra{id\cotensd
f^{-1} \cotensd id} Y\cotensd M,
\end{eqnarray*}
gives the identity on $X\cotensc N$ resp. $Y\cotensd M$, for all
$X\in \comodc$ and $Y\in \comodd$. We have
\begin{eqnarray*}
(id\cotensc id\cotensd f^{-1})(id\cotensc g\cotensc id) &=&
 id\cotensc (id\cotensd f^{-1})(g\cotensc id) \\
 &=& id\cotensc (id\cotensd f^{-1})(id\cotensd f) \\
 &=& id\cotensc (id\cotensd id)\\
 &= & id_{X\cotensc N}.
\end{eqnarray*}
The other composition is similar.
\item is similar to (a).
\end{enumerate}
  \item Similar to (1)
  \item Assume $f$ and $g$ are bijective, let $X\in \comodc$,
  $Y\in \comodd$. We have
\begin{eqnarray*}
FG(X) &= &(X\cotensc N) \cotensd M \\
      &\cong & X\cotensc (N\cotensd M) \\
      &\cong & X\cotensc C \cong X,\quad \text{und} \\
GF(Y) &= & (Y\cotensd M)\cotensc N \\
      &\cong & Y\cotensd (M\cotensc N) \\
      &\cong & Y\cotensd D \\
      &\cong & Y.
\end{eqnarray*}
Hence $F$ is a category equivalence with inverse $G$. Similarly
one can show that $G'$ is an equivalence with inverse $F'$.
\end{enumerate}
\end{proof}
\begin{cor}
Let ($D, C, M, N, f,g$) be a Morita-Takeuchi context.
\begin{enumerate}[(1)]
  \item If $f$ is bijective, then
\begin{enumerate}[(a)]
  \item The comodules $M_C$ and ${}_CN$ are quasi-finite, coflat
  and injectors.
  \item $e_C(M)\cong D$ and $e_C(N)\cong D$ as coalgebras,
  ($e_C(M)$ is the coendomorphism coalgebra of $M$).
  \item $g$ is bijective iff $_CN$ is faithfully coflat iff $M_C$ is
faithfully coflat.
\end{enumerate}
  \item If $g$ is bijective, then
\begin{enumerate}[(a)]
  \item The comodules $N_D$ and ${}_DM$ are quasi-finite, coflat
  and injectors.
  \item $e_D(N)\cong C$ and $e_D(M)\cong C$ as coalgebras.
  \item $f$ is bijective iff $_DM$ is faithfully coflat iff $N_D$ is
faithfully coflat.
\end{enumerate}
 \item If $f$ and $g$ are bijective, then
\begin{enumerate}[(a)]
  \item The comodules $M_C$, $N_D$, $_CN$ and $_DM$ are
  quasi-finite, faithfully coflat and in\-jectors.
  \item $e_C(M)\cong D$, $e_C(N)\cong D$, $e_D(N)\cong C$, and
$e_D(M)\cong C$ as coalgebras.
\end{enumerate}
\end{enumerate}
\end{cor}
\begin{proof}
\begin{enumerate}[(1)]
  \item Assume that $f$ is bijective.
\begin{enumerate}[(a)]
  \item From the fact that $G=-\cotensc N$ is left adjoint to $F=-\cotensd
M$, it follows that $M_C$ is quasi-finite and that $G$ is right
exact. Since $C$ and $N$ are flat $R$-modules we know that $G$ is
also left exact (see \ref{crexact}). So $_CN$ is coflat. Similarly
one can show that $M_C$ is coflat and that $_CN$ is quasi-finite.
The assertion that these comodules are injectors follows from the
fact that a functor between abelian categories which has an exact
left adjoint preserves injective objects.
  \item Follows from the construction of the coendomorphism
  coalgebra.
  \item If $g$ is bijective, then $-\cotensc
N$ is an equivalence, hence $_CN$ is faithfully coflat. If $_CN$
is faithfully coflat, then $g\cotensc id_N : C\cotensc N \lra
(N\cotensd M)\cotensc N$ is an isomorphism, hence $g$ is an
isomorphism. The other assertion is similar.
\end{enumerate}
 \item Similar to (1)
 \item Follows directly from (1) and (2).
\end{enumerate}
\end{proof}
\begin{defn}
A Morita-Takeuchi context ($D, C, M, N, f,g$) in which $f$ and $g$
are both bijective is called \emph{a strict Morita-Takeuchi
context}.
\end{defn}
In lemma \ref{f.inj} we have seen that the involved coalgebras in
a strict Morita-Takeuchi context are Morita-Takeuchi equivalent.
The converse is also true as the following theorem shows
\begin{thm}\label{equ.str}
For two coalgebras $C$, $D$, the following are equivalent
\begin{enumerate}[\emph{(}1\emph{)}]
  \item $C$ is Morita-Takeuchi equivalent to $D$.
  \item There exists a strict Morita-Takeuchi context ($D, C, M, N, f,g$)
\end{enumerate}
\end{thm}
\begin{proof}
We have only to show $(1)\Rightarrow (2)$. Let $F:\comodc
\to\comodd$ be an equivalence with inverse $G:\comodd \to\comodc$,
and let $\phi :id_{\comodc} \lra GF$, $\psi :id_{\comodd}\lra FG$
be the natural transformation belonging to this equivalence. From
\ref{c.equiv} we have $F\cong -\cotensc F(C)$ and $G \cong
-\cotensd G(D)$. Let $f := \phi_C : C\lra F(C)\cotensd G(D)$, and
$g:=\psi_D : D\lra G(D)\cotensc F(C)$.

Claim, ($D, C, G(D), F(C), f, g$) is a strict Morita-Takeuchi
context. The purity and flatness conditions follow from the
general properties of the bicomodules $F(C)$ and $G(D)$. From the
definition of $f, g$ and that the comultiplication of $C$ and $D$
are bicolinear maps it follows that $f$ and $g$ are bicolinear.
The commutativity of the diagrams follows from the relations
\[
F\phi = \psi H, \quad \text{and} \quad G\psi = \phi G.
\]
So we have a strict Morita-Takeuchi context.
\end{proof}
Let $X$ be a right $C$-comodule, and assume that it is
quasi-finite, faithfully coflat and an injector. We will
construct a strict Morita-Takeuchi context derived from $X$. In
\ref{s3} we have seen that $X$ is a $e_C(X)$-$C$-bicomodule and
that $\cohomx{C}$ is a $C$-$e_C(X)$-bicomodule, where
$\cohomx{-}$ is the cohom functor. Now let $f$ be the
$D$-colinear map (see \ref{delta}) $f:=\delta_X : \cohomx{X}\lra
X\cotensc \cohomx{C}$ and $g:= \eta_C : C\lra
\cohomx{C}\cotens{e_C(X)} X$ (see \ref{cunit}).

\begin{prop}\label{scontext}
Let $X$ be as above, then ($e_C(X)$, $C$, ${}_{e_C(X)}X_C$,
$_C\cohomx{C}_{e_C(X)}$, $f$, $g$) is a strict Morita-Takeuchi
context.
\end{prop}
\begin{proof}
It is clear that $X$ is flat in \rmod. From \ref{if.coh.ex} it
follows that $\cohomx{C}$ is also a flat $R$-module, and that
$f=\delta_X$ is bijective. As in the proof of \ref{mtt} it is easy
to show that $g$ is bijective and that $\cohomx{C}\cotens{e_C(X)}
X$ is pure in $\cohomx{C}\tens X$. Since $\cohomx{-}$ is exact,
$X\cotensc \cohomx{C}$ is pure in $X\tens \cohomx{C}$. The
commutativity of the diagrams follows from the defining
properties of $f$ and $g$.
\end{proof}
\begin{cor}
Let $X$ be as above, then $C$ is Morita-Takeuchi equivalent to
$e_C(X)$, the coendomorphism coalgebra of $X$.
\end{cor}
\section{Coalgebras over QF-Rings}\label{s6}
In this section we consider comodule categories for coalgebras
over QF-rings. We will show that Takeuchi's description of the
cohom functor is valid in this case under the condition that the
coalgebra is projective over the ground ring. From now on we
assume that all coalgebras are projective over $R$.

We need the following two lemmas whose proofs can be found in
\cite{Wis98}.
\begin{lem}\label{finite}
Let $D$ be an $R$-coalgebra, $M$ a $D$-comodule. Then
\begin{enumerate}[(1)]
  \item every finite subset of $M$ is contained in a subcomodule
  of $M$ which is finitely generated as $R$-module.
  \item The following are equivalent
\begin{enumerate}
  \item $D$ is finitely generated as $R$-module.
  \item $\comodd = D^{\ast}\textrm{-Mod}$
  ($=\sigma[{}_{D^{\ast}}D]$).
\end{enumerate}
\end{enumerate}
\end{lem}
Let $R$ be a noetherian ring, $D$ an $R$-coalgebra. From lemma
\ref{finite}, it follows that every $D$-comodule is the direct
limit of its subcomodules, that are finitely presented as
$R$-modules. From this and the com-cotensor relations (see
\cite{Altconf}), one obtains some relations between injectivity
and coflatness for coalgebras over QF-rings.
\begin{lem}\label{inj.cof}
Let $D$ be a coalgebra over a QF-ring. A $D$-comodule $M$, which
is flat as $R$-module is injective (resp. an injective
cogenerator) if and only if it is coflat (resp. faithfully coflat)
\end{lem}
From this lemma and \ref{L:homcom} it follows that over a QF-ring
the coalgebra itself is an injective cogenerator as a right and
left comodule.
\begin{cor}\label{injector.inj}
Let $D$ be a coalgebra over a QF-ring, $M$ a right $D$-comodule,
which is flat as $R$-module. If $M$ is injective in $\comodd$,
then it is an injector.
\end{cor}
\begin{proof}
Let $M$ be an injective $D$-comodule, $W$ an injective $R$-module.
Since $R$ is QF, $W$ is flat and $M$ is coflat. Now
\[
(W\tens M)\cotensd - \cong W\tens (M\cotensd -) :\dcomod \to
\rmod.
\]
So the functor $(W\tens X)\cotensd -$ is exact, i.e. $W\tens X$ is
coflat, hence it is injective.
\end{proof}
\begin{prop}\label{end.pres}
Let $D$ be a coalgebra over a coherent ring, $X$ a quasi-finite
$D$-comodule, $P$ a $D$-comodule. If $P$ is finitely presented as
$R$-module, then $\textrm{Com}_D(P,X)$ is finitely presented.
\end{prop}
\begin{proof}
For every index set $\Lambda$, $R^{\Lambda}$ is flat. Consider
\begin{eqnarray*}
R^{\Lambda} \tens \comd(P,X) & \cong & \comd(P,R^{\Lambda}\tens X) \\
     &\cong & \homr(\cohomx{P},R^{\Lambda})\\
     &\cong & \homr(\cohomx{P},R)^{\Lambda}\\
     &\cong & \comd(P,R\tens X)^{\Lambda}  \\
     &\cong & \comd(P,X)^{\Lambda}.
\end{eqnarray*}
the composition of these isomorphisms is the canonical map
\[
R^{\Lambda}\tens \comd(P,X)\to\comd(P,X)^{\Lambda},\; (r_i)\tens f
\mapsto (r_i f)
\]
Hence $\comd(P,X)$ is finitely presented as $R$-module.
\end{proof}
Now we prove our main result in this section, which gives
Takeuchi's representation of the cohom functor for coalgebras over
rings.
\begin{thm}\label{cohom}
Let $C, D$ be coalgebras over a QF-ring $R$, $X$ a
$C$-$D$-bicomodule. If $X_D$ is quasi-finite, injective
cogenerator and $e_D(X)\cong C$ as coalgebras, then the functor
\[
-\cotensc X :\comodc \to\comodd
\]
is an equivalence with inverse
\[
\cohomx{-}\cong -\cotensd \cohomx{D} : \comodd\to \comodc
\]
Moreover, for every $M\in \comodd$, the cohom functor is given by
\[
\cohomx{M}=\lim_{\underset{\Lambda}{\lra}}\textrm{Com}_D
(M_{\lambda},X)^{\ast},
\]
where $\{M_\lambda\}_{\Lambda}$ is the family of subcomodules of
$M$ that are finitely presented as $R$-modules.
\end{thm}
\begin{proof}
From lemma \ref{inj.cof} and corollary \ref{inv.eq} it follows
that $-\cotensc X :\comodc\to\comodd$ is an equivalence with
inverse $\cohomx{-}\cong -\cotensd \cohomx{D} :\comodd\to
\comodc$. For the other assertion consider the following
commutative diagram of categories and functors
\[
\xymatrix{
  \comoddf \ar@<.1cm>[rrr]^{\comd(-,X)}
  \ar@<-.1cm>[rrrdd]_{\cohomx{-}}
    & & & \text{\ccomodf} \ar@<.1cm>[lll]^{\textrm{Com}_C(-,X)}
    \ar@<.1cm>[dd]^{(-)^{\star}} \\
    \\  & & & \comodcf, \ar@<-.1cm>[llluu]_{-\cotensc X}
           \ar@<.1cm>[uu]^{(-)^{\star}}              }
\]
where $\comodcf, \comoddf$ and $\ccomodf$ are the categories of
comodules that are finitely presented as $R$-modules.

Let $P\in \comoddf$, then we have
\begin{eqnarray*}
\comd(P,X)^{\star}\: \cotensc X
    & \cong & \comc(\comd(P,X),X) \\
    & \cong & \comc(X\cotensd P^{\star},X) \\
    & \cong & \comd(P^{\star},\cohomx{D}\cotensc X) \\
    & \cong  &\comd(P^{\star}, D) \cong P,
\end{eqnarray*}
here we have used the facts that $\comd(P,X)$ is finitely
presented, and that $\cohomx{D}\cotensc -$ is right adjoint to
$X\cotensd -$.

For $Q\in \comodcf$ we have
\begin{eqnarray*}
\comd(Q\cotensc X,X)^{\star}
     & \cong & \comc(Q,\cohomx{X})^{\star} \\
     & \cong & \comc(Q,C)^{\star} \\
     & \cong & Q^{\star\star} \cong Q,
\end{eqnarray*}
here we used the fact that $-\cotensc X$ is left adjoint to
$\cohomx{-}$. Hence the functor $\comd(-,X)^{\ast}:\comoddf\to
\comodcf$ is left adjoint to the functor $-\cotensc
X:\comodcf\to\comoddf$ and therefore is isomorphic to
$\cohomx{-}$.

Now let $M\in \comodd$, so
$M=\lim\limits_{\underset{\Lambda}{\lra}}M_{\lambda}$, where
$\{M_{\lambda}\}_{\Lambda}$ is the family of subcomodules of $M$
that are finitely presented as $R$-modules. We have
\begin{eqnarray*}
\cohomx{M} & \cong & \cohomx{\lim_{\underset{\Lambda}{\lra}}M_{\lambda}} \\
           & \cong & \lim_{\underset{\Lambda}{\lra}}\cohomx{M_{\lambda}} \\
           & \cong & \lim_{\underset{\Lambda}{\lra}}
               \comc(M_{\lambda},X)^{\ast}.
\end{eqnarray*}
\end{proof}
Summarizing the results of this section together we get a
characterization for the equivalences of comodule categories over
QF-rings that agrees with Takeuchi's results.
\begin{cor}\label{mteQF}
Let $C, D$ be coalgebras over a QF-ring $R$, $X$ a
$C$-$D$-bicomodule. The following are equivalent
\begin{enumerate}
  \item $X$ is invertible.
  \item The functor $X\cotensd -:\dcomod \to\ccomod$ is an
  equivalence.
  \item $X_D$ is quasi-finite, injective cogenerator and
  $e_D(X)\cong C$ as coalgebras.
  \item ${}_CX$ is quasi-finite, injective cogenerator and
  $e_C(X)\cong D$ as coalgebras.
\end{enumerate}
Moreover, for every $M\in\comodd$, the cohom functor is given by
\[
\cohomx{M}=\lim_{\underset{\Lambda}{\lra}}\textrm{Com}_D
(M_\lambda,X)^{\ast} \cong
\lim_{\underset{\Lambda}{\lra}}(X\cotensd
M_{\lambda}^{\ast})^{\ast},
\]
where $\{M_{\lambda}\}_{\Lambda}$ is the family of subcomodules
of $M$ that are finitely presented as $R$-modules.
\end{cor}
\subsection*{Acknowledgement}
This article is part of my Ph.D thesis at D\"usseldorf
Universit\"at, Germany under the supervision of Prof. Dr. Robert
Wisbauer, many thanks to him for his continuous support and
encouragement. I would also thank Prof. Jos\'e G\'omez-Torrecillas
for the fruitful discussions.
\bibliographystyle{plain}
\bibliography{Takeuchi}

\begin{thebibliography}{10}

\bibitem{Takh99}
K.~Al-Takhman.
\newblock {\em {\"A}quivalenzen zwischen {K}omodulkategorien von {K}oalgebren
  {\"U}ber {R}ingen}.
\newblock Ph. {D}. {D}issertation, H.H.U. D{\"u}sseldorf, Germany, 1999.

\bibitem{Altconf}
K.~Al-Takhman.
\newblock The com and cotensor functors for coalgebras over rings.
\newblock {\em proceedings of the third international palestinian conference on
  mathematics and mathematics education, {B}ethlehem, {P}alestine}, August
  2000.

\bibitem{And73}
F.~Anderson and K.~Fuller.
\newblock {\em Rings and Categories of Modules}.
\newblock Springer-Verlag, Berlin, 1973.

\bibitem{Cae98}
S.~Caenepeel.
\newblock {\em Brauer Groups, {H}opf algebras and {G}alois theory}, volume~4 of
  {\em K-Monographs in Mathematics}.
\newblock Kluwer Academic Publishers, 1998.

\bibitem{Das95}
S.~D{\u{a}}sc{\u{a}}lescu, C.~N{\u{a}}st{\u{a}}sescu, S.~Raianu, and F.~Van
  Oystaeyen.
\newblock Graded coalgebras and {M}orita-{T}akeuchi contexts.
\newblock {\em Tsukuba J. Math.}, 19(2):395--407, 1995.

\bibitem{Gru87}
L.~Grunenfelder and R.~Par\'{e}.
\newblock Families parametrized by coalgebras.
\newblock {\em J. Algebra}, 107:316--375, 1987.

\bibitem{Guzd}
F.~Guzman.
\newblock {\em Cointegration and {R}elative {C}ohomology for {C}omodules}.
\newblock Ph. {D}. {D}issertation, Syracuse University, 1985.

\bibitem{Guz89}
F.~Guzman.
\newblock Cointegration, relative cohomology for comodules, and coseparable
  corings.
\newblock {\em J. Algebra}, 126:211--224, 1989.

\bibitem{Mil65}
J.~W. Milnor and J.~C. Moore.
\newblock On the structure of {H}opf algebras.
\newblock {\em Ann. Math.}, 81:211--264, 1965.

\bibitem{Pop73}
Popesco.
\newblock {\em Abelian categories with applications to rings and modules}.
\newblock Academic Pess, London, 1973.

\bibitem{Sch90}
H.-J. Schneider.
\newblock Principal homogeneous spaces for arbitrary {H}opf algebras.
\newblock {\em Israel J. Math.}, 72(1-2):167--195, 1990.

\bibitem{Sten}
B.~Stenstr{\"o}m.
\newblock {\em Rings and modules of quotients}, volume No. 237 of {\em Lecture
  Notes in Mathematics}.
\newblock Springer Verlag, 1971.

\bibitem{Tak277}
M.~Takeuchi.
\newblock Formal schemes over fields.
\newblock {\em Comm. Algebra}, 5:1483--1528, 1977.

\bibitem{Tak77}
M.~Takeuchi.
\newblock Morita theorems for categories of comodules.
\newblock {\em J. Fac. Sci. Univ. Tokyo}, 24:629--644, 1977.

\bibitem{Wisbi}
R.~Wisbauer.
\newblock {\em Modules and {A}lgebras: {B}imodule {S}tructure and {G}roup
  {A}ctions on {A}lgebras}.
\newblock Pitman Mono. PAM 81, Addison Wesley Longman, Essex, 1996.

\bibitem{Wis98}
R.~Wisbauer.
\newblock Introduction to {C}oalgebras and {C}omodules.
\newblock Lecture Notes, August 1998.

\end{thebibliography}
\end{document}